\documentclass{gtmon_a}
\pdfoutput=1

\proceedingstitle{Groups, homotopy and configuration spaces (Tokyo
  2005)}
\conferencestart{5 July 2005}
\conferenceend{11 July 2005}
\conferencename{Groups, homotopy and configuration spaces,
                in honour of Fred Cohen's 60th birthday}
\conferencelocation{University of Tokyo, Japan}

\editor{Norio Iwase}
\givenname{Norio}
\surname{Iwase}

\editor{Toshitake Kohno}
\givenname{Toshitake}
\surname{Kohno}

\editor{Ran Levi}
\givenname{Ran}
\surname{Levi}

\editor{Dai Tamaki}
\givenname{Dai}
\surname{Tamaki}

\editor{Jie Wu}
\givenname{Jie}
\surname{Wu}

\title{The work of Fred Cohen}

\author{Ran Levi}
\givenname{Ran}
\surname{Levi}
\address{Department of Mathematical Sciences\\University of Aberdeen\\\newline
Meston Building 339\\Aberdeen AB24 3UE\\UK}
\email{ran@maths.abdn.ac.uk}

\volumenumber{13}
\issuenumber{}
\publicationyear{2008}
\papernumber{24}
\startpage{529}
\endpage{546}

\doi{}
\MR{}
\Zbl{}

\arxivreference{}

\keyword{configuration space}
\keyword{braid group}
\keyword{exponents in homotopy theory}
\keyword{mapping class groups}
\keyword{combinatorial group theory}
\keyword{classifying spaces}
\subject{primary}{msc2000}{55R80}
\subject{primary}{msc2000}{55P99}
\subject{secondary}{msc2000}{20F05}
\subject{secondary}{msc2000}{20F36}
\subject{secondary}{msc2000}{55P42}

\received{30 November 2008}
\revised{}
\accepted{5 December 2008}
\published{7 April 2009}
\publishedonline{7 April 2009}
\proposed{}
\seconded{}
\corresponding{}
\version{}

\makeatletter
\def\cnewtheorem#1[#2]#3{\newtheorem{#1}{#3}[section]
\expandafter\let\csname c@#1\endcsname\c@Thm}
\makeatother

\AtBeginDocument{\let\bar\wbar\let\tilde\wtilde\let\hat\what}
\def\longra{\ \hbox{\Large$\longrightarrow$}\ }
\def\rTo{\longra}
\def\rTop#1{\overset{#1}{\longra}}

\newcommand{\F}{{\mathbb F}}
\newcommand{\Fp}{\F_p}   
\newcommand{\pcom}{{}_{p}^{\wedge}}

\newtheorem{Thm}{Theorem}[section]
\cnewtheorem{Prop}[Thm]{Proposition}
\cnewtheorem{Cor}[Thm]{Corollary}
\cnewtheorem{Lem}[Thm]{Lemma}
\cnewtheorem{Claim}[Thm]{Claim}
\cnewtheorem{Ass}[Thm]{Assumption}
\theoremstyle{definition}
\cnewtheorem{Defi}[Thm]{Definition}

\begin{document}
\begin{abstract}
This paper gives an overview of Fred Cohen's work and is a summary of
the talk which I gave during his 60$^{\rm th}$ birthday conference,
held at the University of Tokyo in July 2005.
\end{abstract}

\maketitle

Summarizing Fred's contributions to mathematics in a 1--hour talk, or
for that matter in a single paper, is a daunting task. With nearly a
hundred papers in print, and collaborations with no less than 48
authors, the sheer volume of his publication record makes trying to
choose the appropriate highlights a difficult task. What makes it even
harder is the variety of topics the work touches on.  To all active
topologists Fred Cohen is a very familiar name, and readers are likely
to have a pretty good idea of the main topics his name is mostly
associated with. For instance the study of configuration spaces and
their applications, or his work on homotopy exponents with Moore and
Neisendorfer.  I will of course mention those in this
article. However, as all who know him are aware, Fred's interests are
not limited to mainstream homotopy theory. Quite the contrary, he will
never shy away from an opportunity to explore an unknown, possibly
eccentric grounds.  I will therefore attempt to explore some of the
less known aspects of his work, which I find always interesting, and
occasionally spectacular. This paper is not intended as a
comprehensive summary of Fred's mathematical contributions, but rather
as a sampler of some of his many achievements.

I wish to thank my co-editors in this volume for entrusting me with
the task of delivering a talk on Fred Cohen's work during his 60$^{\rm th}$
birthday conference, and kindly allowing me to record it in print. But
above all I wish to thank Fred himself for simply being there, and
doing so much to enrich our mathematical lives.

\section{Configuration spaces and applications}

Fred's first contributions to the mathematical literature appear in
two research announcements in the Bulletin of the AMS in 1973, as the
outcomes of his PhD thesis \cite{C1,C2}.  The real work however only
appeared three years and five independent publications later in the
collection by Cohen--Lada--May \cite{CLM}, which became more or less
immediately the definitive basic reference in the subject of its title
-- the homology of iterated loop spaces.  To this colossal 490 page
book, Fred contributed two articles with a total page length of 192
pages.

Very few papers in algebraic topology are more fundamental than these
two. The following is the tip of an iceberg summary of what is done
there.

{\it ``The homology of $\mathcal{C}_{n+1}$--spaces, $n\geq 0$'' },
\cite{C3}, starts with a careful analysis of the Dyer--Lashof and Browder
homology operations, the relationships between them and the
Pontrijagin product in the homology of iterated loop space. The main
idea is the utilization of May's little $n$--cubes operad and its
action on the homology of an iterated loop space. This is then applied
to the calculation of the homology of May's configuration space model
$H_*(C_{n+1}X)$, and $H_*(\Omega^{n+1}\Sigma^{n+1}X, \Fp)$, as free
objects, in the appropriate sense, on the homology of $X$. I can't
possibly explain these results here in much detail, nonetheless how
they are obtained. But to give the reader the flavour, here is a brief
description.

\begin{Defi} A graded $\Fp$--vector space $M$ is:
\begin{itemize}
\item an allowable $R_n$--module if
there are homomorphisms
\[Q^s \co M_q \to M_{q+2s(p-1)}, \quad (Q^s\co M_q\to
M_{q+s}\quad\mbox{if p=2})\] for $0\le 2s< q+n$, ($s<q+n$) such that $Q^s=0$ for
$2s<q$,  $(s<q)$ and the composition of the $Q^s$ satisfy the Adem
relations.

\item an allowable $AR_n$--module if, in addition $M$ admits an action of
the dual of the Steenrod algebra which satisfies the Nishida
relations.

\item an allowable $AR_n$--algebra, if in addition $M$ is a commutative
algebra satisfying \begin{enumerate} \item  $Q^sx=x^p$, if $|x|=2s$, $(|x|=s)$ for any
$x\in M$, \item $Q^s(1) = 0$ if $s>0$ and \item the Cartan formula for
products. \end{enumerate}

\item an allowable $AR_n\Lambda_n$--Hopf algebra (with conjugation), if $M$ is a
monoidal Hopf algebra satisfying further properties concerning
the Browder operations $\lambda_n$ and the ``top'' operation $\xi_n$, and
their commutation relations with the Dyer Lashof operations and the
Pontrijagin product.
\end{itemize}
\end{Defi}

 Define
\begin{itemize} 
\item[(1)]a functor $W_n$ to be the free functor left adjoint to the
forgetful functor from allowable $AR_n\Lambda_n$--Hopf algebras to
cocommutative coalgebras over the dual of the Steenrod algebra.

\item[(2)] A functor $G$ from allowable $AR$--Hopf algebras to $AR$--Hopf
algebras with conjuation  to be
\[GW = W\otimes_{\pi W} G_0W,\]
where $W$ is an allowable $AR$--Hopf algebra, $\pi W$ is the
commutative monoid under the product in $W$, $\pi GW$ is the
commutative group generated by $\pi W$, and $G_0W$ is its group
ring. In other words $GW$ is the localization of the ring $W$ at the
monoid $\pi W$.
\end{itemize}

With this terminology and notation, Fred proves the following:

\begin{Thm}{\rm\cite[Theorem 3.1]{C3}}\qua For every space $X$ there is an
isomorphism of allowable $AR_n\Lambda_n$--Hopf algebras:
\[\bar\eta_*\co W_nH_*(X)\rTo H_*(C_{n+1}X)\]
\end{Thm}

\begin{Thm}{\rm\cite[Theorem 3.2]{C3}}\qua For every space $X$ there is an
isomorphism of allowable $AR_n\Lambda_n$--Hopf algebras with
conjugation:
\[\tilde\eta_*\co GW_nH_*(X)\rTo H_*(\Omega^{n+1}\Sigma^{n+1}X)\]
\end{Thm}

The generality of the results implies of course that certain specific
cases, particularly when one is not after full and complete answers,
may be easier to compute directly than by referring to the
theorems. But nevertheless, the importance of these results in the
development of algebraic topology cannot be over estimated.

In {\it ``The Homology of $SF(n+1)$''}, \cite{C4}, Fred studies the
homology of the topological monoid of degree 1 self maps of the sphere
$S^{n+1}$, where the monoid operation is given by composition. The
main theorem is the statement that the Pontrijagin ring $H_*(SF(n+1),
\F_p)$ is a commutative algebra which for odd primes is isomorphic as
an algebra to $H_*(\Omega^{n+1}_0S^{n+1})$ (where the loop space
structure is given as usual by juxtaposition of loops. The proof is by
direct calculation, and utilizes similar techniques to those used in
\cite{C3}.

\section{Braid groups}

Artin's braid groups are fundamental objects which arise naturally in
geometry, knot theory, ring theory and many other mathematical
disciplines. They also have prominent roles to play in algebraic
topology, to a large extent thanks to Fred's contributions.  Fred's
first published encounter with Braid groups appear as a 2 page
appendix to \cite{C3}, where he computes the homology of the braid
group on $r$ strands $B_r$ as a module over the dual of the Steenrod
algebra, as well as the rational and integral homology.  The mod--$2$
results read as follows.

\begin{Thm} {\rm\cite[Theorem A.1]{C3}}\qua There is an isomorphism of modules over
the dual of the Steenrod algebra
\[H_*(B_r,\F_2) \cong P[\xi_{2^j-1}\,|\,j\ge 1]/I,\]
where $I$ is the ideal generated by the monomials
\[\xi_{j_1}^{k_1}\cdots\xi_{j_t}^{k_t}\quad\mathrm{where}\quad
\sum_{i=1}^tk_i2^{j_i} > r.\]
Furthermore, the action of the dual Steenrod algebra is determined
by the requirements that $Sq_*^r$ acts trivially if $r>1$, and that
$Sq_*^1(\xi_{j+1}) = \xi_j^2$. 
\end{Thm}

Mark Mahowald has shown in a 1977 Topology paper \cite{Mah} that the
Thom spectrum of the natural map $\bar{\eta}\co\Omega^2S^3 \rTo BO$
is the Eilenberg--MacLane spectrum $K(\Z/2,0)$. It is common in
mathematics to look for different ways in which various bits of
mathematical knowledge fit together. Fred has always been on the
lookout for such 'non-accidents'. In \cite{C5} he shows the following:

\begin{Thm}{\rm\cite{C5}}\qua Let $B_\infty$ denote the colimit of the braid groups
$B_r$ under the obvious inclusions. Then there is a homology
isomorphism $\theta\co K(B_\infty, 1) \rTo \Omega^2S^3$.
\end{Thm}

From this he is able to deduce the following beautiful

\begin{Cor}{\rm\cite{C5}}\qua The Thom spectrum $MB_\infty$ of the composite
$\bar{\eta}\circ \theta$ (at any prime) is the $K(\Z/2,0)$--spectrum.
Thus every mod--2 homology class may be realized as a manifold whose
stable normal bundle has a $B_r$--reduction.
\end{Cor}

Fred's romance with the braid groups continued through his career and
to the present day. One beautiful result was presented by him in this
conference, based on the the paper \cite{CW} in this volume.

Another very recent work of a totally different flavour is a
collaboration with Alejandro Adem and Dan Cohen \cite{ACC}. A group
$\Gamma$ is said to be homologically toroidal if there is a
homomorphism $\Z^{n_1}*\cdots *\Z^{n_k}\rTo \Gamma$ inducing a split
epimorphism on integral homology. An example of such a group is the
pure braid group $P_r$.

The authors show the following:

\begin{Thm}{\rm\cite{ACC}}\qua If $\Gamma$ is a homologically toroidal group 
and $\Gamma\to U(n)$
is a representation then the composite
\[B\Gamma\rTo BU(n)\rTo BU\]
is null homotopic, while if $\Gamma\rTo O(n)$ is a representation
the corresponding composite is null homotopic if and only if the
first two Stiefel Whitney classes of the representation vanish.
\end{Thm}

In the process, the subgroup of elements in the $K$--theory of
$B\Gamma$ which arise from orthogonal representations is also
determined.

For the pure braid group, each quadratic relation in the cohomology
ring $H^*(P_r)$ is shown to correspond to a spin representation of
$P_r$. This representation is nontrivial, but it gives rise to a
trivial bundle over the configuration space $F(\C,r) = K(P_r, 1)$.

\section{Exponents in homotopy theory}
The mid to late 70s saw a collaboration between Fred, Joe Neisendorfer
and John Moore, the results of which mark some of the most beautiful
results in unstable homotopy theory ever achieved. The subject of
study is exponents in homotopy theory.  Homotopy theorists, realizing
that looking for explicit calculations in unstable homotopy groups
couldn't possibly be feasible in great generality, started looking for
qualitative, rather than quantitative results. The search for exponent
results was one of the paths one could explore. Two famous
conjectures, both open to this day, are worth mentioning, as they
provided much of the motivation for the Cohen--Moore--Neisendorfer
project.

The Barratt conjecture states that if a double suspension $X=\Sigma^2
Y$ has the property that the order of the class of the identity
element in the abelian group $[X,X]$ is $p^r$ for some prime $p$, then
the $p^{r+1}$ power map on $\Omega^2 X$ is null homotopic. The Moore
conjecture is more general in its setup, but less specific in its
conclusion. It states that if $X$ is a finite $p$--local $CW$ complex,
then the torsion part of $\pi_*(X)$ has a global exponent if and only
if the rational homotopy of $X$ is globally finite dimensional, or
using Moore's terminology, if and only if $X$ is elliptic (as opposed
to hyperbolic). Both conjectures were stated at a point were not a
single example was known. A lucid discussion of these conjectures can
be found in \cite{Sel}.

The Cohen--Moore--Neisendorfer team set out in the mid 70's to fix the
situation.  The three of them together published five papers on the
subject, two of which appeared in the Annals of Mathematics
\cite{CMN1,CMN2,CMN3,CMN4,CMN5}.

In the 1983 International Congress of Mathematicians in Warsaw, Fred
delivered an invited address where he reported on the
Cohen--Moore--Neisendorfer project.  His report appeared as \cite{C6}.

The Cohen--Moore--Neisendorfer papers contain enough ideas to keep a
whole generation of topologists busy. They influenced the work of many
topologists in the almost 30 years since the first paper was
published. Among them Anick, Gray, Theriault, and Selick.  Exponents
in homotopy theory were studied before this project commenced, but
arguably never before in such a systematic fashion. One of the most
striking aspect of the project is the elegant and systematic use of
techniques of differential graded Lie algebras. The authors apply
these methods to the homotopy Bockstein spectral sequence, which is a
differential graded Lie algebra with respect to the Samelson
product. Using the homological information, they conclude the
existence of a product splittings of certain loop spaces. The first
theorem we quote is an example. For a prime $p$ and a positive integer
$m$, let $P^n(p^m)$ denote the homotopy cofibre of the degree $p^m$
map on the sphere $S^{n-1}$, and let $S^{n-1}\{p^m\}$ denote the
homotopy fibre of the same map In \cite{CMN1,CMN2} the authors
restrict attention to primes $p>3$ and prove the following.

\begin{Thm}{\rm\cite{CMN1}}\qua Let $p$ be an odd prime, and $n$ a positive integer. Then
\[ \Omega P^{2n+2}(p^r) \simeq S^{2n+1}\{p^r\}\times
\Omega\left(\bigvee_{m=0}^\infty P^{4n+2mn+3}(p^r)\right).\]
\end{Thm}

Another important ingredient in their analysis is the homotopy fibre
$F_n\{p^r\}$ of the pinch map $P^n(p^r)\rTo S^n$. In two consecutive
theorems in \cite{CMN1} they provide a product splitting for $\Omega
F_n\{p^r\}$.

\begin{Thm}{\rm\cite{CMN1}}\qua Let $p$ be an odd prime, and $n$ a positive integer. Then \[\Omega F^{2n+1}\{p^r\} \simeq S^{2n-1} \times
\prod_{k=1}^\infty S^{2p^kn-1}\{p^{r+1}\} \times P,\] and
\[\Omega F^{2n}\{p^r\} \simeq \Omega S^{2n-1}\times S^{4n-3} \times
\prod_{k=1}^\infty S^{4p^kn-2p^k-1}\{p^r\}\times P'.\] Here $P$ and
$P'$ denote (different) infinite products of loop spaces on mod $p^r$
Moore spaces.
\end{Thm}

These product splittings allows the authors to construct a map
\[\pi \co \Omega^2
S^{2n+1} \rTo S^{2n-1},\] whose composition with the double
Freudenthal suspension map
\[S^{2n-1}\rTop{E^2} \Omega^2 S^{2n+1}\rTop{\pi} S^{2n-1}\]
is homotopic to the degree $p$ map on $S^{2n-1}$.  Combined with work
of Toda this allows them to show that the $p$--torsion in
$\pi_*(S^{2n+1})$ has exponent $p^{n+1}$.

In the same paper they also show that the $\pi_*(P^n(p^r))$ contains
infinitely many elements of order $p^{r+1}$.

A refinement of the their methods in \cite{CMN2} allows them to chose
a map $\pi$ as above, such that the composition the other way
\[\Omega^2S^{2n+1} \rTop{\pi} S^{2n-1} \rTop{E^2} \Omega^2
S^{2n+1}.\] is the double loops of the $p^{\rm th}$ power map on
$\Omega^2S^{2n+1}$.

All that is needed now is to iterate this composite $n$ times, and the
$p^n$ power map on $\Omega^2S^{2n+1}$ factors through $S^1$. This
implies that $p^n$ annihilates the homotopy of $S^{2n+1}$.

In the next two papers \cite{CMN3} and \cite{CMN4}, they show the
existence of exponents for the Moore spaces $P^n(p^r)$. In \cite{CMN3}
they prove that $p^{2r+1}$ annihilates the homotopy of an even
dimensional mod $p^r$ Moore space at odd primes. To improve on this,
as to fit with the Barratt conjecture which predicts an exponent
$p^{r+1}$, the missing ingredient is a product splitting for $\Omega
P^{2n+1}(p^r)$. This is done in \cite{CMN4}, and in a subsequent paper
by Neisendorfer, the predicted exponent is obtained.

Examples of the Moore and Barratt conjectures have been constructed by
numerous authors following Cohen--Moore--Neisendorfer, but as general
statements these conjectures remain as intact today as they were when
they were originally stated.

\section{An early curiosity}

In the 70's between his thesis and the work he did with Moore and
Neisendorfer, Fred wrote a number of papers where he applied his
configuration spaces techniques to various problems -- most notably a
powerful generalization of the Borsuk--Ulam theorem. The contributions
already mentioned above are among Fred's greatest achievements. In
this short section however, I have chosen to mention a much less
familiar paper he collaborated on early in his career.

Every mathematician who has been active long enough knows what Erd\"os
numbers are. Paul Erd\"os wrote according to MathSciNet more than 1500
papers, most of them in collaboration with other mathematicians. This
uncommon prolificacy yielded the concept. A mathematician has Erd\"os
number 1 if they wrote a paper with Erd\"os himself, and Erd\"os
number $\le n$ if they collaborated with a mathematician whose Erd\"os
number is $n-1$. It is conjectured that any mathematician who ever
collaborated on a paper has a finite Erd\"os number. However, since
Erd\"os was a number theorist, it appears unlikely for a homotopy
theorist to have a small Erd\"os number.  Fred, never in the habit of
putting a title on what interested him, earned an Erd\"os number 2 due
to a collboration with number theorist Selfridge in \cite{CS}.  The
main theorem in this paper is of the kind Erd\"os himself would
probably approve of.

\begin{Thm}{\rm\cite{CS}}\qua
There exist infinitely many odd numbers $M$, such that neither $M+2^n$
not $|M-2^n|$ is a prime power for any $n$.
\end{Thm}

The authors also construct an explicit 94 digit example of the
theorem.

\section{``You'd think it's easy to decide whether something is
divisible by 2''}

The title of this section is actually a quote of Fred in his 1990
algebraic topology class, which I was a part of. This is of course all
about the so called \emph{strong form of the Kervaire invariant
conjecture}, a statement which mystified and deceived topologists for
more than 50 years.  The statement is simple:
\medskip

{\it The Whitehead square
$\omega_{2n+1}=[\iota_{2n+1},\iota_{2n+1}]\in\pi_{4n+1}(S^{2n+1})$ is
divisible by 2 if $n=2^k-1$.}
\medskip

It is easy to see $\omega_{2n+1}=0$ for $n=0, 1, 3$ (since the spheres
		     in question are $H$--spaces).  In \cite{C9}, which
		     is a textbook treatment of this and many other
		     aspects of classical homotopy theory, Fred gives
		     no less than 5 equivalent formulations of this
		     question.  In the next statement cohomology is
		     taken with coefficients in $\F_2$.
\begin{Thm}{\rm\cite[Proposition 11.4]{C9}}\qua
Let $n\neq 0,1,3$. Then the following statements are equivalent.
\begin{enumerate}
\item The Whitehead product $\omega_{2n+1}$ is divisible by 2.
\item The short exact sequence
\[0 \rTo \Z/2 \rTo \pi_{4n+l}(S^{2n+1})\rTo \pi_{4n+2}(S^{2n+2})\rTo 0\]
is not split .
\item There is a map $P^{4n+2}(2) \rTo \Omega S^{2n+2}$ which is
  non-zero in homology
\item There exists a space X with $\bar{H}^i(X) = \Z/2$ for $i=2n+2,
  4n+3$, and $4n+4$, and zero otherwise, with $Sq^{2n+2}\co  H^{2n+2}(X)
  \rTo H^{4n+4}(X)$ and $sq^l\co  H^{4n+3}(X)\rTo H^{4n+4}(X)$
  isomorphisms.
\item $\Omega^2[-1]$ is homotopic to $-1$ on $\Omega^2S^{2n+l}$.
\end{enumerate}
\end{Thm}

For graduate students in Rochester Fred's special ``homework problems''
are a familiar concept. I remember vividly how after proving the
theorem above in class he suggested, ``well, here is a homework
problem for you. Prove one of these statements''. Anybody who ever
heard a lecture by Fred is likely to have been assigned a homework
problem in the context of the Kervaire invariant conjecture, and many
other subjects. Those who know him are aware of three basic facts: (1)
He is genuinely interested to know the answer, (2) he has tried it
himself, and (3) you may spend your lifetime trying to solve this
homework problem, and there is no partial credit.

But Fred tends to be very serious about his homework assignments, and
when the students are struggling he always tries to help. So if five
formulations are not enough, then in \cite{C10} he gives yet another
formulation. This one is quite special in that it relates the question
to the real Cayley--Dickson algebras. He constructs a certain subspace
$K(n,\epsilon)$ of the topological vector space given by the
polynomial ring $\R[x,y]$. He shows the following:

\begin{Lem}{\rm\cite{C10}}\qua If $\epsilon>0$ and $n\ge 2$, the space
$K(n,\epsilon)$ is homotopy equivalent to the $(4n-1)$--skeleton of
$\Omega^2S^{2n+1}$.
\end{Lem}

He then uses the multiplication induced from the Cayley--Dickson
algebra to construct a model $Sq$ for the degree 2 map on $S^k$. He
then proceeds to show

\begin{Lem}{\rm\cite{C10}}\qua If $k=2^n-1$ and $n\geq 3$, then $\omega_k$ is
divisible by 2 if the loop squaring map and the map $\Omega^2(Sq)$ are
homotopic when restricted to the $(4n-3)$--skeleton.
\end{Lem}

Using these two lemmas, Fred concludes the following

\begin{Prop}{\rm\cite{C10}}\qua If $k=2^n-1$, then $\omega_k$ is divisible by 2
if and only if the composites
\[K(\frac{k-1}{2},\epsilon)\rTo
\Omega^2S^k\rTop{\ \Omega^2(Sq)\ }\Omega^2S^k,\quad\mathrm{and}\]
\[K(\frac{k-1}{2},\epsilon)\rTo
\Omega^2S^k\rTop{2\ }\Omega^2S^k\] are homotopic.
\end{Prop}
This gives a way (which unfortunately fails) to attempt an explicit
homotopy that does the job.

\section{Some general homology calculations}
Fred's work includes many computational results, some specific, and
others very general. It is the second kind this section deals
with. One example of such a calculation is of course Fred's thesis,
where he gives a complete description of the homology of
$\Omega^n\Sigma^nX$, but there are many others, some of which I will
touch on below.

Many problems in algebraic topology involve understanding mapping
spaces. The most obvious example is that of iterated loop spaces. The
$n$--fold loop space $\Omega^n X$ can be identified with the pointed
mapping space $\mathrm{Map}_*(S^n, X)$. Along totally different lines,
the Sullivan conjecture, and subsequent work by Miller and Lannes
involve studying mapping spaces of the form $\mathrm{Map}(BV, X)$,
where $V$ is an elementary abelian $p$--group.

When $X$ and $Y$ are arbitrary spaces, identification of
$\mathrm{Map}_*(X,Y)$ is practically impossible. In a joint work with
Larry Taylor \cite{CT}, the authors study these spaces under certain
hypotheses on $X$ and $Y$, for which they obtain a rather explicit
result.

\begin{Thm}{\rm\cite{CT}}\qua Let $Y$ be an $m$--fold suspension, and let $X$ be a
finite complex of dimension less than $m/2$, which is itself the
suspension of a connected space. Then there is a mod--$p$ homology
isomorphism of graded vector spaces \[H_*(\mathrm{Map}(X,Y))\approx
\bigotimes H_*((\Omega^i Y)^{\beta_i(X)}),\] where $\beta_i(X)$ is the
$i^{\rm th}$ Betti number of $X$, and the tensor product runs over all $i$
such that $\beta_i(X)\neq 0$. Furthermore, if $X$ is also a double
suspension, then the isomorphism is as Hopf algebras.
\end{Thm}

Another family of spaces which features frequently in Fred's work are
configuration spaces $C^k(M)$, already mentioned in the context of his
thesis.  Fred's interest in configuration spaces never withered, and
they keep coming up in his work in a variety of contexts.  Here are a
few examples.

In \cite{BCT} the authors B\"odigheimer, Taylor and Fred study the
configuration space $F(M,k)$ of $k$ distinct points in a smooth
compact $m$--manifold $M$, possibly with boundary. The paper determines
the additive structure of the homology $H_*(F(M,k);\F)$ where $\F$ is
any field if $m$ is odd, and $\F_2$ otherwise. This is a well cited
paper, but like many important results its significance was not
discovered immediately. This paper appeared in 1989, was first cited
in 2000, and since then twelve more times, none of which by any of the
authors themsleves.

In \cite{CT2} the authors study the cohomology of the configuration
space $F(\R^m,r)$ as a module over the symmetric group
$\Sigma_r$. Although most of the work is done with integer
coefficients, the most specific results, including identification of
specific characters, are obtained with rational coefficients.  The
authors identify a class $A_{2,1}\in H^{m-1}(F(\R^m, 2);\Z)$. Letting
$\pi_{i,j} \co  F(\R^m, r) \rTo F(\R^m, 2)$ be defined by the formula
\[\pi_{i,j}(x_1, \ldots , x_r) = (x_i, x_j),\] they define $A_{i,j} =
\pi_{i,j}(A_{2,1})$. These classes have many good properties and they
play a key role in of $H^*(F(\R^m, r);\Z)$. This paper is quoted by a
number of authors in various applications, notably in a recent paper
by Arone, Lambrechts and Volic \cite{ALV}.

Finally, in a collaboration with Sam Gitler \cite{CG}, the loop space
homology of $F(M,k)$ for certain manifolds $M$ is studied. Most of the
work is concerned with the case where $M$ is obtained by removing a
single point from a closed manifold. For instance, if $M = \R^m$, $m
\ge 3$, the authors prove that the primitive elements of the integral
homology ring $S = H_*(\Omega F(\R^m, k))$ form a Lie algebra
generated by elements $B_{i,j}$ subject to the infinitesmal braid
relations and that $S$ itself is the universal enveloping algebra of
this Lie algebra. The classes $B_{i,j}$ are related to the generators
used by Fred in the calculation of $H^*(F(\R^m, k))$ in his
thesis. More generally, if $M$ is a simply connected punctured
manifold of dimension $m\ge 3$, then the authors show that $\Omega
F(M, k)$ is homotopy equivalent to a direct product, one factor of
which is $\Omega F(\R^m, k)$. The remaining factors, which are
identified explicitly and involve $M$, are also loop spaces, but the
equivalence in the product decomposition is not multiplicative. Thus,
to describe the ring structure in homology, one must determine the
twisting among the factors, for instance, between the various
$B_{i,j}$ and classes coming from the homology of $M$. Under
appropriate hypotheses on $M$ and with coefficients in certain fields
$\F$, the authors find a complete set of relations describing this
twisting, thereby determining $H^*(\Omega F(M, k); \F)$
completely. They also show that the hypotheses are necessary by giving
examples where some of the relations do not hold.

\section{Mapping class groups}

The homotopy theory associated to the mapping class groups is another
subject Fred contributed very substantially to, with no less than nine
papers to his name with the phrase ``mapping class group'' in the
title. In these papers Fred and his coauthors explore various
connections of the mapping class groups to homotopy theory, or perform
various cohomology calculations (a good example of the latter is
\cite{BeC}).

A very beautiful example of the way Fred explores connections among
mathematical objects is in his paper \cite{C11}. Let $M^g$ be a closed
orientable surface of genus $g$ and let $\Gamma_g$ denote its mapping
class group. The hyperelliptic mapping class group $\Delta_g$ is
defined to be the centralizer in $\Gamma_g$ of the hyperelliptic
involution which acts on $M^g$ and fixes $2g +2$ points. In this paper
Fred studies the groups $\Delta_g$. For $g=2$, $\Delta_d=\Gamma_g$,
but for $g>2$ these subgroups are neither normal nor of finite index
in $\Gamma_g$.  Let $\Gamma^n$ denote the mapping class group of
$S^2$, with $n$ fixed points. The group $\Gamma^n$ was studied from
the group theoretic point of view by Magnus. The relevance to
$\Delta_g$ comes from the existence of a central extension
\[0\rTo\Z/2 \rTo \Delta_g \rTo \Gamma^{2g+2}\rTo 1.\]
In this paper Fred uses techniques of classical homotopy theory to
study topological and homological properties of $\Delta_g$. In
particular he constructs spaces of type $K(\pi, 1)$, where $\pi
\Delta_g$. The constructions involves properties of the Lie groups
$SO(3)$ and $Spin^c(3)$, and particularly a work of S. Smale, who
showed that the natural inclusion $SO(3)\rTo \mathrm{Diff}^+(S^2)$ is
a homotopy equivalence.

\section{Combinatorial group theory in homotopy theory}\label{simp}

Fred's work on combinatorial group theory in homotopy theory is, in my
mind, one of his most beautiful and original contributions.  The core
paper, for reasons he must know better than I do, remains unpublished
\cite{C12}. When I say, the work remains unpublished, I'm lying a bit.

The object of study in this paper is the group $[\Omega\Sigma
X,\Omega\Sigma X]$ of pointed homotopy classes of maps, where $X$ is
any reasonable space. Within this group one can single out two types
of elements:

\begin{enumerate}
\item For each natural number $k$ the class of the $k^{\rm th}$ power map
 \[\Omega\Sigma X \rTop{\Delta} (\Omega\Sigma X)^{\times k}\rTop{\mu}
 \Omega\Sigma X,\] where $\Delta$ is the $k$--fold diagonal map, and
 $\mu$ is the loop space multiplication map.
\item For each natural number $k$ the class of the composite
\[\Omega\Sigma X\rTop{h_k}\Omega\Sigma
X^{(k)}\rTop{\Omega\omega_k}\Omega\Sigma X,\] where $h_k$ is the
$k^{\rm th}$ James--Hopf invariant, $X^{(k)}$ is the $k$--fold smash power of
$X$, and $\omega_k$ is the $k$--fold iterated Whitehead product.
\end{enumerate}

Naturality of these maps implies that they can be considered as
endomorphisms (ie, self natural transformations) of the functor
$\Omega\Sigma$ on the homotopy category of spaces. The set of all such
endomorphisms forms a group under loop multiplication at the target
(rather than under composition, which gives a monoid, but not
generally a group, structure).  In this group, which we denote
$[\Omega\Sigma(-), \Omega\Sigma(-)]$, one can consider the subgroup
$H_\infty$ generated by the elements above.

Let $J_n(X)$ denote the $n^{\rm th}$ stage of the James construction on
$X$. Thus, $J_\infty(X) = \cup_n J_n(X)$ is homotopy equivalent to
$\Omega\Sigma X$. Then, as above, one can consider the group (under
loop multiplication) of homotopy classes of natural transformations
from $J_n(-)$ to $\Omega\Sigma(-)$, which we denote by $[J_n(-),
\Omega\Sigma(-)]$. The James construction is naturally filtered by
subfunctors $J_n(-)$, and this filtration induces a filtration on the
group $H_\infty$ by subgroups $H_n\le H_{n+1}\le\cdots\le
H_\infty$. As Whitehead products, Hopf invariants, and compositions of
such are among the most important maps in classical homotopy theory,
studying the group $H_\infty$ may give a breakthrough in our
understanding of these maps and the relationships between them.

A similar exercise can be done in the group of natural transformations
in the homotopy category $[(-)^n, \Omega\Sigma(-)]$ (again, under loop
multiplication). In this group consider the subgroup generated by the
elements $p_i$, given as the classes of the composites
\[(-)^n\rTop{{\rm proj}_i} (-)\rTop{E} \Omega\Sigma(-).\]
Let $K_n\le [(-)^n, \Omega\Sigma(-)]$ be the subgroup generated by
these.

Here is one of the fundamental theorems proven in \cite{C12}.

\begin{Thm}{\rm\cite{C12}}\qua The group $K_n$ is a finitely presented, torsion free
nilpotent group of class $n$. A specific presentation is given by the
group generated by elements $x_1, \ldots x_n$, subject to the
relations:
\begin{enumerate}
\item[\rm(i)] $[x_{i_1},\ldots, x_{i_k}] = 1$ if $x_{i_j}=x_{i_k}$ for
$j<k$.
\item[\rm(ii)] $[x_{i_1}^{n_1}, \ldots, x_{i_k}^{n_k}] = [x_{i_1},\ldots,
x_{i_k}]^{n_1\cdots n_k}$.
\end{enumerate}
\end{Thm}

He then proceeds by identifying the groups $H_n$ as subgroups given by
a given set of generators inside the groups $K_n$.

Although \cite{C12} is not published, it has inspired further works by
J Wu and others, some in collaboration with Fred.

\section{Classifying spaces  -- a personal note}\label{personal}

Classifying spaces are not a main theme in Fred's work, but like so
many other subjects it is one on which he touched and inspired others
-- in this case myself. A mathematician is measured by his work and
contributions to mathematics, but the hard work of teaching and
inspiring students is often neglected. To do so in Fred Cohen's case
will be to miss out on what I think is a major aspect of his
mathematical persona. Of course, I only have my own experiences as
Fred's student to share, but I dare guess that my story is not
atypical.

In 1989 I went to graduate school in Rochester. My advisor at the
time, Emmanuel Farjoun, recommended it to me very highly as one of the
best places in the world to do a PhD in Topology. He told me who was
there, and gave me a brief description of each person and his work. On
Fred he said that he's been through some hard times health-wise, but
in spite of that he was a wonderful mathematician.

I got to Rochester, and I remember very vividly the very first class
with Fred, whom until then I never met. We were all seated waiting for
him to arrive.  Based on Emmanuel's description, my own prejudice made
me expect a weak and tortured figure.  Instead, in came a man,
anything but weak and tortured, almost running with the aid of his
cane, with this huge smile on his face, and a strong confident voice.
Said hello, and started one of the most illuminating lectures I've
ever heard until then. It was at this point, I think, that I decided
this man will be my thesis advisor.

During my first year in Rochester I had ample opportunity to talk
mathematics to Fred. I loved his lectures. A condition of
participation in Fred's classes was that each one of us had to solve
at least one homework problem in public, and he gave us plenty to
think about. His manner was deceptively very casual. I remember
thinking -- oh what a wonderful new way of proving statements by saying
``well, what could it be?'' in a convincing tone. It didn't quite occur
to me at that stage that one has to be prepared to explain why it
couldn't be anything else. The first time I tried this technique on
Fred in a private session, he replied ``Ran, it could be many
things. Go back and work out exactly what it is.''. And, of course, he
had a good reason for that. I had confidently ``proved'' a very wrong
statement. This was very educational. Fred is always very friendly and
informal with students, but they have to get it right or meet his
``steel''.

My own work with Fred was an example of how wide his interests
range. In 1990 he participated in the Barcelona conference, which took
place in the beautiful beach town of San Feliu on the Catalan Costa
Brava. Trying to study a question of Dave Benson which occurred in a
discussion between them in the Hotel's bar, Fred ``gave birth'' to my
thesis subject.

Two background concepts before we explain the question. A group is
said to be ``perfect'' if its first integral homology vanishes. If $G$
is a perfect group and $X$ is a space with $\pi_1(X)=G$, then the
Quillen's ``plus'' construction associates with $X$ a simply connected
space $X^+$ with the same homology as $X$.

Benson's question was: what can be said on the homotopy type of
$BG^+$, where $G$ is a finite perfect group? Fred reacted by doing one
of the things he does best. He looked at a few examples of finite
perfect groups to which he applied $B(-)^+$, took the corresponding
loop spaces, and calculated the living daylight out of them. Within
the course of the evening, he managed to calculate a few examples
which exhibited quite a curious behavior. They were all ``finitely
resolvable by fibrations over spheres and loop spaces on spheres''. In
other words he produced a finite sequence of fibrations, the total
space in the first of which is $\Omega BG^+$, where the fibre in the
$n^{\rm th}$ fibration is the total space in the $(n+1)^{\rm st}$, and
where all base spaces where either spheres or loop spaces on
spheres. He recorded his thoughts in a little paper, which he
published in the Conference Proceedings \cite{C7}.

Here is an easy example. Let $p$ be an odd prime, and let $n\ge 2$ be
an integer dividing $p-1$. Then the cyclic group $\Z/n$ act on $\Z/p$
by automorphisms and one can form the semidirect product $G(p,n) =
\Z/p\rtimes \Z/n$. This group is not perfect, but it is $p$--perfect,
and one can replace the ``plus'' construction by $p$--completion for a
$p$--local version. The observation is that
\[\Omega BG(p,n)\pcom \simeq S^{2n-1}\{p\},\]
where the right hand side is the fibre of the degree $p$ map on the
sphere.  Thus one has a length 2 resolution
\[\Omega S^{2n-1} \rTo \Omega BG\pcom\rTo S^{2n-1}.\]
This is of course a very easy example. The article \cite{C7} contained
quite a few more, some of them far from obvious. This was convincing
enough for Fred to make what he called a ``Rush Conjecture'', that in
general $\Omega BG^+$ is finitely spherically resolvable.

When he came back to Rochester and told me about this amusing
discovery, I was totally fascinated -- in fact, more than fascinated, I
was hooked. I always loved group theory, especially finite, and the
chance to work on a combination of group theory and homotopy theory
seemed too good to let pass. So I almost immediately asked him to
become my advisor and to let me work on this question.

Since I spent many hours before this moment in time with Fred
discussing mathematics, I was very surprised, not to mention
disappointed, when he wasn't fast to agree.  Not only he wasn't going
to let me work on this wonderfully eccentric subject; he didn't want
to be my advisor, or so I thought. I was crushed! A day later he left
a note in my mailbox telling me that he would in fact be happy to talk
to me about whatever I wanted, including becoming my thesis advisor,
if I insisted. Later I learnt the reason for his initial reluctance.
Fred always regarded being an advisor as a great responsibility,
almost a type of fatherhood. He was genuinely concerned, and for a
good reason, about things like getting a job after graduation, and
about certain problems being too bizarre for a PhD project, and this
one was certainly an example of such a problem. So, in a sense he
maybe wanted to make sure that his own students know exactly what it
is they are getting themselves into. It seems to me that in most if
not all cases, they knew. I certainly did.

Fred's ``Rush Conjecture'' became my thesis subject, and I kept at it
for a number of years after graduating. It turned out to be a subject
much richer and more interesting than could have been predicted during
that pub chat with Benson, and the following Proceedings article. The
conjecture itself turned actually to be wrong, as I proved about a
year after graduation. There are examples, in fact rather easy
examples, of finite $p$--perfect groups which are not spherically
resolvable \cite{L}.  However other aspects of these spaces remain
very interesting, and inspired a number of other mathematicians. One
remarkable example is a recent work of Benson \cite{Be}, where he
gives a purely algebraic interpretation of the mod $p$ loop space
homology of $BG\pcom$.

Fred and I wrote a few more papers together, on classifying spaces and
other subjects. Several other of his students have shared the same
pleasure with me. These were and still are illuminating and fruitful
interactions, for which I am thoroughly grateful. At the time of
writing this summary it is too late to wish Fred a happy birthday, but
I will conclude by saying:

\begin{center}
{\Large All the best to you Fred,

for many years of mathematics to come.}
\end{center}

\bibliographystyle{gtart} \bibliography{link}

\end{document}